# Distorted stochastic dominance: a generalized family of stochastic orders

Tommaso Lando[1,2*] and Lucio Bertoli-Barsotti[1]

[1] University of Bergamo, Department of Management, Economics and Quantitative Methods, Via dei Caniana 2, Bergamo, Italy.

[2] VŠB-Technical University of Ostrava, Department of Finance, Sokolskà Trida 33, Ostrava, Czech Republic.

*Corresponding author

E-mail: tommaso.lando@unibg.it

**Abstract.** We study a generalized family of stochastic orders, semiparametrized by a distortion function $H$, namely $H$-distorted stochastic dominance, which may determine a continuum of dominance relations from the first- to the second-order stochastic dominance (and beyond). Such a family is especially suitable for representing a decision maker's preferences in terms of risk aversion and may be used in those situations in which a strong order does not have enough discriminative power, whilst a weaker one is poorly representative of some classes of decision makers. In particular, we focus on the class of power distortion functions, yielding power-distorted stochastic dominance, which seems to be particularly appealing owing to its computational simplicity and some interesting statistical interpretations. Finally, we characterize distorted stochastic dominance in terms of distortion functions yielding isotonic classes of distorted expectations.

**Funding.** This research was supported by the Czech Science Foundation (GACR) under project 17-23411Y (to T.L.).

1. Introduction

In the theory of decision under uncertainty, decision makers measure their preferences regarding uncertain prospects by assigning different weights, to be interpreted either as misjudgements or as subjective revisions, to the outcomes of the corresponding random variable (RV) or to the corresponding probabilities. Mathematically, this weighting process may be formulated as a transformation of the values of the RV or of its cumulative distribution function (CDF), which, in turn, may be expressed, for instance, in terms of *integrals* (integrated CDFs, integrated quantiles, etc.), *utilities* (functions of the RV) or *probability distortions* (functions of the CDF). Based on different combinations of such transformations, the theory of stochastic dominance (SD) provides tools for representing the preferences of decision makers and their attitudes towards risk.

In this context, the most commonly used SD relations are *first-* and *second-order stochastic dominance* (FSD, SSD, respectively), owing to their several applications in areas such as economics, econometrics, finance and insurance. Basically, FSD represents any decision maker who prefers "more" to "less", whereas SSD represents any decision maker who is also "risk averse". Although most decision makers may be represented by FSD, it is generally not easy to establish whether one uncertain prospect is "bigger" than another (checking FSD is generally a strong condition). Thus, the discriminative power of FDS is generally poor. On the other hand, SSD might be limiting for those decision makers who are "mostly" risk averse but may have some degree of flexibility in their preferences and therefore exhibit a weak risk attitude, at least in some situations (as discussed by Muller et al. 2017). This may be illustrated by a paradoxical example: consider a choice between 1) a sure gain of $a$ and 2) an uncertain gain of 0 (with probability 0.5) or $a + \varepsilon$ (with probability 0.5), in which we assume $\varepsilon$ to be positive and arbitrarily small compared with $a$. Clearly, all risk-averse decision makers would choose option 1 (for small $\varepsilon$s). Nevertheless, we argue that most of the others, who are usually not risk averse, would make the same choice. Hence, neither FSD nor SSD represents those who choose option 1. Proper justifications for decisions of this type may be provided by a general approach, making it possible to generalize both FSD and SSD as well as yielding SD relations "between" these two.

The literature contains various examples of SD relations that interpolate FSD and SSD. Fishburn (1976, 1980) established continua of SD relations for bounded and unbounded probability distributions that fill the "gap" between FSD, SSD and weaker SD relations of integer degrees. Leshno and Levi (2002) defined the *almost SD* of the first and second order, which allows for small violations of the FSD or SSD rules, whereby the weight of such violations is controlled by a real parameter. More recently, some authors focused on this topic again. Tzeng et al. (2013) proposed an adjustment for the main theorem in the paper by Leshno and Levi (2002). Tsetlin et al. (2015) generalized second-order almost SD to dominance rules of a higher degree. Muller et al. (2017) introduced a new and different family of stochastic orders, covering preferences from FSD to SSD.

Drawing inspiration from these works, we are also concerned with finding orders between FSD and SSD, but we look at the problem from a totally different perspective. In fact, the various approaches of the papers cited above are all related to the concept of a utility function, that is, a transformation, or weighting, of the outcome of an RV. Differently, we refer to the dual approach of Yaari (1987), which focuses on a transformation of the corresponding CDF, referred to as probability distortion. We recall that a distortion is an increasing function $H: [0,1] \to [0,1]$ such that $H(0) = 0$ and $H(1) = 1$. This transformation process may be seen as a decision maker performing a subjective weighting of the original CDF, in which the choice of $H$ may represent different ways of measuring uncertainty (the decision maker's perceptions). For instance, a concave distortion function emphasizes the weight of smaller outcomes, which conforms to the idea of risk aversion, whereas a convex $H$ emphasizes the weight of the larger ones.

Levi and Wiener (1998) studied the SD relations between distorted distributions and investigated the classes of distortions that preserve FSD and SSD. Following their approach, in section 2 we compare RVs of which the CDFs are transformed through a common distortion function $H$ and consequently define a semiparametric family of stochastic orders, denoted by *H-distorted stochastic dominance (H-DSD)*, $\geq_H$. We study the relationships among the orders $\geq_H$ for different choices of $H$. Basically, the strength of $H$-DSD is related to the shape of the distortion function and especially to its degree of concavity/convexity. If $H$ is "more convex" than $G$ in the sense of Chan et al. (1990), then $H$-DSD implies $G$-DSD. By varying the degree of convexity of $H$, we may establish a continuum of SD relations, from FSD to SSD and beyond (i.e., weaker than) SSD. This can be achieved by focusing on a parametric family of distortion functions $H_k(t) = H(t, k)$ and, consequently, by defining a parametric family of stochastic orders. In particular, we choose the class of power functions $H_k(t) = t^k, k > 0$, which gives rise to the *power-DSD (PDSD)* of order $1 + 1/k$, where the order determines the strength of the SD relation and $1/k \in (0, \infty)$ determines its degree of risk aversion. In section 3, we show that PDSD satisfies some desirable conditions, yielding FSD as a limiting case and SSD as a special case.

Insofar as $H$-DSD generalizes SSD, in section 4, we apply the same approach to generalize the *increasing and convex order* (ICX), an order that is somewhat complementary to SSD (Shaked and Shanthikumar, 2007). Differently from SSD, the ICX represents decision makers who prefer "more" to "less" but are also "risk lovers". We define a generalization of ICX, via a distortion function $H$, and denote it as $H$-*risk-loving-DSD*. Given two (possibly) different distortions $H_1, H_2$, $H_1$-DSD and $H_2$-risk-loving-DSD can be combined to define the $H_1, H_2$-*mixed-DSD* order. Similarly to what has been undertaken recently by Muller et al. (2017), $H_1, H_2$-mixed-DSD imposes constraints on aversion as well as attraction to risk, expressed in terms of $H_1, H_2$.

In section 5, we characterize DSD, risk-loving-DSD and mixed-DSD in terms of classes of order-preserving functionals, generally known in the literature as distorted expectations or distortion risk measures (see for instance Wang and Young 1998). The generalized Gini indices, introduced by Donaldson and Weymark (1983), also belong to this class. Wang and Young (1998) showed that distorted expectations derived from convex distortion functions preserve ICX. Similarly, it can be shown that distorted expectations derived from concave distortion functions preserve SSD. We finally show that a distorted expectation preserves $H$-DSD ($H$-risk-loving-DSD) if it derives from a distortion that is less (more) convex than $H$, whereas mixed-DSD is preserved if the distortion function satisfies some intuitive constraints in terms of convexity and concavity.

## 2. Distorted stochastic dominance

Let $f$ be a non-decreasing and right continuous function. We define the right-continuous (generalized) inverse of $f$ as $f^{-1}(p) = \sup\{z: f(z) \leq p\}$ (Marshall et al. 2011, p. 714). We aim to compare a pair of RVs $X$ and $Y$ with corresponding CDFs $F_X$ and $F_Y$ and quantile functions $Q_X = F_X^{-1}$ and $Q_Y = F_Y^{-1}$. To avoid some technical issues, in this paper, we consider only RVs with finite expectations. We recall that the definitions of FSD and SSD can be expressed equivalently in terms of the quantile function. In fact, FSD and SSD are equivalent to first- and second-degree *inverse stochastic dominance* (Thistle 1989).

**Definition 1.** We say that $X$ FSD dominates $Y$ and write $X \geq_1 Y$ iff
$$F_X(t) \leq F_Y(t), \forall t \in \mathbb{R}$$
or, equivalently,
$$Q_X(t) \geq Q_Y(t), \forall t \in [0,1].$$

**Definition 2.** We say that $X$ SSD dominates $Y$ and write $X \geq_2 Y$ iff

$$\int_{-\infty}^{u} F_X(t)dt \leq \int_{-\infty}^{u} F_Y(t)dt, \forall u \in \mathbb{R}$$

or, equivalently,

$$\int_0^p Q_X(t)dt \geq \int_0^p Q_Y(t)dt, \forall p \in [0,1].$$

A new semiparametric family of stochastic orders may be obtained by comparing via SSD pairs of RVs of which the CDFs are transformed by a common distortion function $H$. The corresponding order is denoted as $H$-DSD. The definition can be stated with some alternative formulations, based on the equivalence results of Proposition 1 below.

For an RV, $X$, with CDF $F_X$, we denote with $X_H$ the RV with CDF $F_{X_H} = H(F_X)$.

**Proposition 1**

The following conditions are equivalent.

1) $X_H \geq_2 Y_H$.

2) $\int_0^u Q_X(H^{-1}(t))dt \geq \int_0^u Q_Y(H^{-1}(t))dt, \forall u \in [0,1]$.

3) $\int_0^u Q_X(t)dH(t) \geq \int_0^u Q_Y(t)dH(t), \forall u \in [0,1]$.

4) $\int_{-\infty}^x H(F_X(t))dt \leq \int_{-\infty}^x H(F_Y(t))dt, \forall x \in \mathbb{R}$.

Now we can provide the definition of $H$-DSD.

**Definition 3.** We say that $X$ dominates $Y$ w.r.t. $H$-DSD and write $X \geq_H Y$ iff any of the equivalent conditions of Proposition 1 holds true.

$H$-DSD is equivalent to SSD between $H$-distorted distributions. Moreover, it can be seen that $H$-DSD is closely related to (weak) p-majorization (Marshall et al. 2011, p. 583) or, using a different terminology, to the *weak spectral order* (Chong 1974), with respect to the measure associated with $H$.

The distortion process, studied, among others, by Levi and Weiner (1998), Wang and Young (1998) and Yaari (1987), may be interpreted as a subjective weighting of the probabilities, which reflects the attitude (e.g. towards risk) of decision makers. Intuitively, a concave (convex) $H$ emphasizes left (right) tail probabilities and, indirectly, attaches more weight to the corresponding smaller (larger) outcomes, conforming to the idea of risk aversion (attraction). For instance, a risk-averse decision maker whose probability of a negative event (loss, failure, etc.) is given by $p$ will basically act as if such a probability was greater than $p$.

The idea of obtaining a family of stochastic orders via distortion functions has been studied elsewhere. Li and Shaked (2007) made use of distortion functions to generalize the *total time on test transform order*. Levi and Wiener (1998) proved that all distortion functions preserve FSD, whilst all and only concave distortion functions preserve SSD. Put otherwise, $X \geq_1 Y$ implies $X_H \geq_1 Y_H$ (this result is obvious), whilst $X \geq_2 Y$ implies $X_H \geq_2 Y_H$ iff $H$ is concave (Levi and Wiener 1998, Theorem 6).

The next theorem establishes the relation between the orders $\geq_H$ for different $H$s. We show that the strength of $H$-DSD is related to the degree of concavity/convexity of the distortion function or, equivalently, to its Arrow–Pratt risk aversion measure (Arrow 1971; Pratt 1964). We provide some useful definitions below. Definition 4 formalizes the concept of one function being more convex (or concave) than another (Chan et al. 1990; Van Zwet 1964). Definition 5 gives the expression of the Arrow–Pratt measure.

**Definition 4 (Chan et al. 1990).** Let $F, G$ be a pair of CDFs. We say that $F$ is more (less) convex than $G$ and write $F \geq_{CX} G$ ($G \geq_{CX} F$) iff $F(G^{-1})$ is convex (concave) in the interval $[0,1]$.

**Definition 5.** Define the Arrow–Pratt risk aversion measure of a function $f$, at least twice differentiable, as

$$R_f(u) = -\frac{f''(u)}{f'(u)}.$$

**Theorem 1**

Let $H, G$ be two distortion functions, at least twice differentiable.

1) The following conditions are equivalent.

   i)   $G'(u)/H'(u)$ is decreasing in $[0,1]$.
   ii)  $R_G(u) \geq R_H(u), \forall u \in [0,1]$.
   iii) $H \geq_{CX} G$.

2) If any of the conditions above (i, ii or iii) holds, $X \geq_H Y$ implies $X \geq_G Y$.

For $H(t) = t$, $H$-DSD is equivalent to SSD. In particular, the identity function represents a watershed for the DSD family in that concave (convex) distortions yield orders that are weaker (stronger) than SSD. Differently, it can readily be seen that $X \geq_1 Y$ implies $X \geq_H Y$, for all distortion functions $H$, but, conversely, no distortion function exists, say $\tau$, such that $X \geq_\tau Y$ iff $X \geq_1 Y$. Hence, FSD represents an upper bound for the DSD family. Intuitively, for Theorem 1, as the degree of convexity of $H$ increases, we might come "close" to this upper bound. Then, our idea is to focus on parametric families of distortion functions $H_k(t) = H(t, k)$. In particular, we search for families that fulfil the following conditions.

*C1) Antisymmetry.* $X \geq_{H_k} Y$ and $Y \geq_{H_k} X$ implies $F_X = F_Y$.

*C2) Monotonicity.* If $k_1 > k_2$, $X \geq_{H_{k_1}} Y$ implies $X \geq_{H_{k_2}} Y$.

*C3) Identity.* For some value of $k$, $H_k(t) = t$.

*C4) Consistency.* A set K exists such that $X \geq_{H_k} Y$, $\forall k \in K$, iff $X \geq_1 Y$ (equivalence with FSD).

C1 implies that $H_k$ must assign positive weight to all the probabilities in the interval $[0,1]$ to satisfy the antisymmetry property, also denoted as "neutrality" by Yaari (1987). According to Theorem 1, C2 implies that parameter $k$ must determine the degree of concavity/convexity of $H_k$ and consequently the strength of the dominance relation. C3 states that the family $H_k$ must contain the identity function to generalize SSD. C4 states that, if $H_k$-DSD holds for all $k \in K$, then it must be equivalent to FSD ($H_k$-DSD implies FSD, whereas the converse implication is always true).

In section 3 we show that the class of power distortion functions satisfies all the above conditions.

### 3. Power-DSD

We search for a family of distortion functions that fulfils C1–C4. It is not difficult to identify classes of distortions that satisfy C1–C3. Thus, C4 is crucial. We show that the class of power distortion functions satisfies them all, besides being very simple and providing useful interpretations.

Let $H_k(t) = t^k$ ($k > 0$). $H_k$ yields a parametric family of stochastic orders, which we denote as PDSD of order $1 + 1/k$.

**Definition 6.** We say that $X$ dominates $Y$ w.r.t. the PDSD of order $1 + \frac{1}{k}$ and write $X \geq_{\left(1+\frac{1}{k}\right)-PD} Y$ iff $X \geq_{H_k} Y$, where $H_k(t) = t^k$ ($k > 0$). In particular, we write $X \geq_{1-PD} Y$ iff $X \geq_{H_k} Y, \forall k > 0$.

It can easily be seen that PDSD satisfies C1–C3. Below we prove that it also satisfies C4.

**Lemma 1**

Let $H_k(t) = t^k$ and let $Q$ be a quantile function. Then

$$\lim_{k \to \infty} \frac{\int_0^u Q(t) dH_k(t)}{H_k(u)} = Q(u)$$

holds at every point of continuity of $Q$.

Lemma 1 implies the following result.

**Theorem 2**

$X \geq_{1-PD} Y$ iff $X \geq_1 Y$.

Thus, PDSD fulfils C1–C4; in particular, we obtain:

C1) $X \geq_{\left(1+\frac{1}{k}\right)-PD} Y$ and $Y \geq_{\left(1+\frac{1}{k}\right)-PD} X$ iff $F_X = F_Y$.

C2) If $k_1 > k_2$, $\geq_{\left(1+\frac{1}{k_1}\right)-PD}$ implies $\geq_{\left(1+\frac{1}{k_2}\right)-PD}$.

C3) $X \geq_{2-PD} Y$ iff $X \geq_2 Y$.

C4) $X \geq_{1-PD} Y$ iff $X \geq_1 Y$.

PDSD also has an interesting statistical interpretation. Assume that $k$ is a positive integer and denote with $U_{k:k}$ the largest-order statistic (or sample maximum) from a sample of i.i.d. RVs $U_1, \ldots, U_k$. The CDF of $U_{k:k}$ is $F_{U_{k:k}}(t) = F_U(t)^k$. Clearly, we obtain

$$X \geq_{\left(1+\frac{1}{k}\right)-PD} Y \text{ iff } X_{k:k} \geq_2 Y_{k:k}, \tag{1}$$

which, in turn, implies $X_{h:h} \geq_2 Y_{h:h}$ for all $h \leq k$ (C1). As $k$ increases, larger values become progressively more important (whereas the weights of smaller values remain fixed) and PDSD approaches FSD.

For $k > 1$, PDSD covers the preferences of decision makers from FSD to SSD, similarly to the family of stochastic orders recently introduced by Muller et al. (2017). Nevertheless, for $k < 1$, we obtain orders that are weaker than SSD. Such weaker orders can be used to increase the rate of completeness, to be understood as the proportion of pairs of distributions that are ranked according to a given preorder, providing finer criteria for the decision making (see for instance Muliere and Scarsini 1989). In a decision problem, $k$ may be chosen according to the degree of risk aversion of the decision maker, which may be quantified by the number $1/k \in (0, \infty)$. In this case, $1/k$ might be seen as a latent parameter to be determined or inferred from data, for instance through economic experiments or surveys with questions on hypothetical gambles (Anderson and Mellor 2009). Conversely, given two RVs $X, Y$, we may be interested in finding the strongest PDSD relation between $X$ and $Y$, that is, the largest $k$, say $k^* = k^*(X, Y)$, such that $X \geq_{\left(1+\frac{1}{k^*}\right)-PD} Y$. Knowledge of the exact value of $k^*$ provides information about the actual strength of the dominance relation between $X$ and $Y$ (the higher, the stronger). In section 3.1, we describe how to determine $k^*$ in some special cases. It is not always possible to have an explicit solution, in terms of $k$, for the PDSD inequality condition, but, in any case, $k^*$ can be found numerically. In particular, PDSD is simple to verify in the discrete case (as we show in section 3.1) and for single-crossing distributions, owing to Theorem 3 below.

We say that $F_X, F_Y$ are single-crossing (from below) if a point, say $x_1$, exists such that $F_X(x) \geq F_Y(x)$ for $x > x_1$ and $F_X(x) \leq F_Y(x)$ for $x < x_1$. In this special case, PDSD verification reduces to a comparison of the generalized expectations $E(X_{H_k})$, $E(Y_{H_k})$, where $E(U_{H_k}) = \int_0^1 Q_U(t) dt^k$ ($H_k(t) = t^k$). This can be stated as follows (the result follows from Theorem 3 of Hanoch and Levy 1969).

**Theorem 3**

If $F_X, F_Y$ are single-crossing (from below), then $X \geq_{\left(1+\frac{1}{k}\right)-PD} Y$ iff $E(X_{H_k}) \geq E(Y_{H_k})$ ($H_k(t) = t^k$).

The statistical interpretation of Theorem 3 is as follows. If $k$ is a positive integer and $F_X, F_Y$ are single-crossing (from below), $X \geq_{\left(1+\frac{1}{k}\right)-PD} Y$ iff $E(X_{k:k}) \geq E(Y_{k:k})$. Then, $X$ dominates $Y$ if $F_X$ starts below $F_Y$ and the sample maximum (of order $k$) of $X$ is expected to be larger than that of $Y$.

### 3.1. Some examples of PDSD

In the following, we apply PDSD to both discrete (e.g. empirical case) and continuous RVs.

For discrete distributions, the computation of $\int_0^u Q_X(p)dH_k(p)$ is quite simple. Let $X$ be a discrete RV that takes values $x_1, \ldots, x_n$ (e.g., a set of empirical observations). The QF of $X$ is

$$Q_X(p) = x_i \text{ for } p \in I_i, i = 1, \ldots, n,$$

where $I_i = [F_{i-1}, F_i) = [F_X(x_{i-1}), F_X(x_i))$ and $F_0 = 0$. Then,

$$\int_0^u Q_X(p)dH_k(p) = \sum_{i=1}^{j-1} x_i ((F_i)^k - (F_{i-1})^k) + x_j \left(u^k - (F_{j-1})^k\right), \text{ for } u \in I_j.$$

(From this representation, it is easy to verify Lemma 1: $\lim_{k\to\infty} u^{-k} \int_0^u Q_X(p)dp^k = \sum_{i=1}^n x_i \mathbf{1}_{I_i}(p) = Q_X(p)$, almost everywhere — the limit moves to zero in the "jump" points of $Q_X$, that is, for $u = F_{i-1}, i = 1, \ldots, n$).

As for continuous distributions, we focus on location scale families and use the single-crossing argument of Theorem 3.

*Example 1.* A particularly interesting case of single-crossing distributions is the comparison of an RV $X$ with a constant, $c$. If $c$ is contained in the support of $X$, then the CDFs of $X$ and $c$ cross once from below; thus, the RV $X$ dominates the certain value $c$, hence $X \geq_{\left(1+\frac{1}{k}\right)-PD} c$, iff $E(X_{H_k}) \geq c$, that is, if the expectation of the distorted RV $X_{H_k}$ is greater than (or equal to) $c$.

Let $X = c$, where $0 \leq c \leq 1$ and $Y$ are a uniform RV defined on the support $[0,1]$. Clearly, $X$ is less "risky" than $Y$, but the two CDFs always cross so that $X \not\geq_1 Y$ (except for the case $c = 1$). If, for instance, $c = 0.5$, then $X \geq_2 Y$, so we argue that, for $c > 0.5$, the dominance relation between $X$ and $Y$ is stronger than SSD, whereas, for $c < 0.5$, the dominance relation between $X$ and $Y$ is weaker. $E(X_{H_k}) \geq c$ for $\frac{k}{k+1} \geq c$, hence $X \geq_{1/c-PD} Y$; that is, the order of the PDSD is the reciprocal of $c$. This result confirms our conjecture. For instance, if $c = 0.2$, $X \geq_{5-PD} Y$, if $c = 0.6$, $X \geq_{1.\overline{6}-PD} Y$ and, if $c = 0.8$, $X \geq_{5-PD} Y$. Clearly, if $c = 1$, then $X \geq_1 Y$. We also study the behaviour for the $(1+\gamma)$-SD order defined by Muller et al. (2017), namely $\geq_{(1+\gamma)-SD}$. We obtain that $X \geq_{(1+\gamma)-SD} Y$ for $\gamma \geq \frac{(1-c)^2}{c^2}$, where $\gamma \in [0,1]$; thus, for instance, if $c = 0.6$, $X \geq_{1.\overline{4}-SD} Y$ and, if $c = 0.8$, $X \geq_{1.0625-SD} Y$.

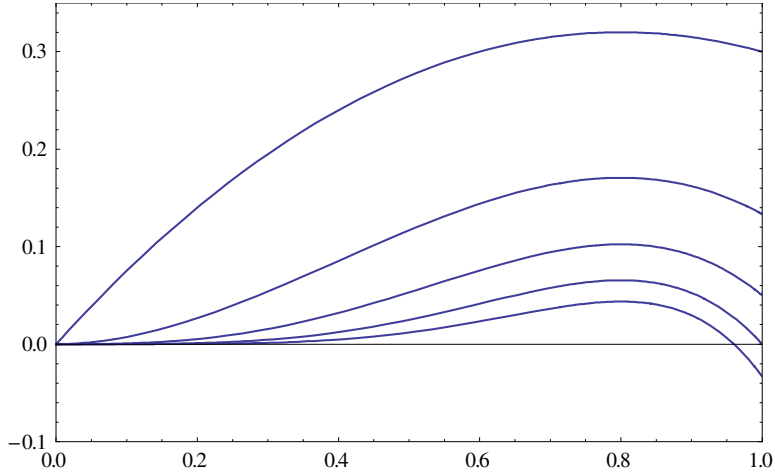

**Figure 1.** $\int_0^p Q_X(t)dt^k - \int_0^p Q_Y(t)dt^k$ for $c = 0.8$ and $k = 1,\ldots,5$. For higher values of $k$, the curve becomes lower. The function is always non-negative (i.e., PDSD holds) for $k \leq 4$.

*Example 2.* Let, for instance, $X$ take the values 1,3,6 (with probabilities $\frac{1}{2}, \frac{1}{4}, \frac{1}{4}$) and $Y$ take the values 0,2,4,5 (with uniform probability $\frac{1}{4}$). The two CDFs cross twice and $E(X) = E(Y) = 2.75$. It is easy to verify that $X \not\geq_2 Y$; thus, we may search for a dominance relation that is weaker than SSD. We can easily compute

$$\int_0^u Q_X(p) - Q_Y(p)dp^k = \begin{cases} t^k & 0 < t \leq 1/4 \\ 2^{1-2k} - t^k & 1/4 < t \leq 3/4 \\ t^k - 2^{1-2k}(3^k - 1) & 3/4 < t \leq 1 \end{cases}$$

and find $\int_0^u Q_X(p) - Q_Y(p)dp^k \geq 0, u \in [0,1]$, for $k \leq \frac{\ln 2}{\ln 3} = 0.63$. Thus, $X \geq_{2.58-PD} Y$.

*Example 3. Location–scale families.* Let $X, Y$ belong to the same location–scale family. Then $Q_X(t) = G^{-1}(t)\sigma_X + \mu_X$ and $Q_Y(t) = G^{-1}(t)\sigma_Y + \mu_Y$ where $G$ is the CDF of an RV, say $Z$, with location parameter 0 and scale parameter 1. If $\mu_X \geq \mu_Y$ and $\sigma_X < \sigma_Y$, then $F_X, F_Y$ are single-crossing from below. The crossing point is at $p = G\left(\frac{\mu_X - \mu_Y}{\sigma_Y - \sigma_X}\right)$. Thus, for Theorem 3, $X \geq_{\left(1+\frac{1}{k}\right)-PD} Y$ iff

$$E_H(X) = E_H(Z)\sigma_X + \mu_X \geq E_H(Z)\sigma_Y + \mu_Y = E_H(X)$$

If $\mu_X = \mu_Y$, $X \geq_H Y$ iff $\sigma_X < \sigma_Y$ and $E_H(Z) \leq 0$. If $\mu_X \neq \mu_Y$, $X \geq_{\left(1+\frac{1}{k}\right)-PD} Y$ iff

$$E_H(Z) \leq \frac{\mu_X - \mu_Y}{\sigma_Y - \sigma_X}. \tag{2}$$

*Logistic distribution*

Let $X \sim Log(\mu_X, \sigma_X)$ and $Y \sim Log(\mu_Y, \sigma)$. $S^-(F_X - F_Y) = 1$ (with sign sequence $-, +$) for $\mu_X \geq \mu_Y, \sigma_X \leq \sigma_Y$. Then, (2) reduces to

$$E_{H_k}(Z) = \int_0^1 \ln\left(-1 + \frac{1}{x}\right) kp^{k-1} dp = \gamma + \Psi(k) \leq \frac{\mu_X - \mu_Y}{\sigma_Y - \sigma_X},$$

where $\gamma$ is the Euler's constant and $\Psi(k) = (\ln \Gamma(k))'$ is the digamma function.

The inequality can be solved numerically. Let $\mu_X = 0.1, \mu_Y = 0, \sigma_X = 1, \sigma_Y = 1.1$. $E_{H_k}(Z) \leq 1$ for $k \leq 2$; thus, we find that $X \geq_{1.5-PD} Y$. Clearly, we find identical results for $\mu_X = 1, \mu_Y = 0, \sigma_X = 1, \sigma_Y = 2$. Meanwhile, for $\mu_X = 1, \mu_Y = 0, \sigma_X = 1, \sigma_Y = 1.1$ (or, identically, $\mu_X = 0.1, \mu_Y = 0, \sigma_X = 1, \sigma_Y = 2$), $E_{H_k}(Z) \leq 1$ for $k \leq 12367.46$; thus, we find that $X \geq_{1.00008-PD} Y$, which is very close to FSD (actually, the two QFs cross at $1/(1 + e^{1/10}) \cong 1$). Put otherwise, the sample maximum of $X$ is expected to be larger than that of $Y$, for random samples up to 12367. Now, let $\mu_X = \mu_Y = 0, \sigma_X = 1, \sigma_Y = 1.1$. It is known that, in this case, $X \geq_2 Y$. In particular, SSD is the strongest dominance relation between $X$ and $Y$. In fact, $E_{H_k}(Z) \leq 0$ for $k \in [0,1]$; thus, we find that $X \geq_{\left(1+\frac{1}{k}\right)-PD} Y$ holds $\forall k \in [0,1]$.

*Normal distribution*

Let $X \sim N(\mu_X, \sigma_X)$ and $Y \sim N(\mu_Y, \sigma)$. $S^-(F_X - F_Y) = 1$ (with sign sequence $-, +$) for $\mu_X \geq \mu_Y, \sigma_X \leq \sigma_Y$. Then, (2) becomes

$$E_{H_k}(Z) = \int_0^1 \text{erf}^{-1}(2p) kp^{k-1} dp \leq \frac{\mu_X - \mu_Y}{\sigma_Y - \sigma_X}.$$

The inequality can be solved numerically. Let $\mu_X = 0.1, \mu_Y = 0, \sigma_X = 1, \sigma_Y = 1.1$. $E_{H_k}(Z) \leq 1$ for $k \leq 3.82$; thus, we find that $X \geq_{1.26-PD} Y$. Interestingly, we obtain identical results for $\mu_X = 1, \mu_Y = 0, \sigma_X = 1, \sigma_Y = 2$. Whilst, for $\mu_X = 1, \mu_Y = 0, \sigma_X = 1, \sigma_Y = 1.1$ (or, identically, $\mu_X = 0.1, \mu_Y = 0, \sigma_X = 1, \sigma_Y = 2$), the $\geq_{\left(1+\frac{1}{k}\right)-PD}$ approaches FSD (actually, the two QFs cross at $1/2 \, \text{erf}(-5\sqrt{2}) \cong 1$). Now, let $\mu_X = \mu_Y = 0, \sigma_X = 1, \sigma_Y = 1.1$. It is known that, in this case, $X \geq_2 Y$. In particular, SSD is the strongest dominance relation between $X$ and $Y$. In fact, $E_{H_k}(Z) \leq 0$ for $k \in [0,1]$; thus, we find that $X \geq_{\left(1+\frac{1}{k}\right)-PD} Y$ holds $\forall k \in [0,1]$.

## 4. Risk-loving- and mixed-DSD

As shown in section 3, $H$-DSD is equivalent to SSD between distorted distributions; that is, $X \geq_H Y$ iff $X_H \geq_2 Y_H$. This approach can be extended to other SD relations, besides SSD. In particular, we focus on the increasing convex order, an order that is somewhat complementary to SSD, usually referred to as the *increasing convex order* (*ICX*).

**Definition 7.** We say that $X$ dominates $Y$ w.r.t. ICX and write $X \geq^2 Y$ iff

$$\int_u^\infty 1 - F_X(t)dt \leq \int_u^\infty 1 - F_Y(t)dt, \forall u \in \mathbb{R}$$

or, equivalently,

$$\int_p^1 Q_X(t)dt \geq \int_p^1 Q_Y(t)dt, \forall p \in [0,1]$$

or, equivalently, $-Y \geq_2 -X$.

SSD and ICX differ just for the verse of integration. By performing an integration from the right, ICX is intended to attach more weight to the right tail of the distribution and less to the left one, which clearly represents attraction towards risk. For this reason, we may denote ICX as *risk-loving-SSD*. Correspondingly, we can define a *risk-loving* distorted order w.r.t. a distortion function $H$.

**Definition 8.** We say that $X$ dominates $Y$ w.r.t. *risk-loving-$H$-DSD* and write $X \geq^H Y$ iff $X_H \geq^2 Y_H$. Equivalently, $X \geq^H Y$ iff $-Y \geq_H -X$.

Let $\widetilde{H}(t) = 1 - H(1-t)$ be the *dual distortion function* of $H$. Note that, integrating by substitution,

$$\int_0^p Q_{-X}(u)dH(u) = \int_0^p -Q_X(1-u)H'(u)du = \int_{1-p}^1 Q_X(y)\widetilde{H}'(y)du = \int_{1-p}^1 Q_X(y)d\widetilde{H}(y).$$

Thus, $\geq^H$ can also be defined equivalently as follows:

$$X \geq^H Y \text{ iff } X_{\widetilde{H}} \geq^2 Y_{\widetilde{H}}. \tag{3}$$

Clearly, the properties of risk-loving-DSD are closely related to those of DSD. In particular, Theorem 1 also holds for $\geq^H$; that is:

$$\text{If } H \geq_{CX} G, X \geq^H Y \text{ implies } X \geq^G Y. \tag{4}$$

$\geq^H$ and $\geq_H$ have the same basic properties. For $H(t) = t$, $\geq^H$ is equivalent to $\geq_H$. Moreover, similarly to what we discussed for $\geq_H$, $\geq^H$ with convex (concave) $H$ is stronger (weaker) than ICX

($\geq^2$). Owing to the relation among $\geq^H$ and $\geq_H$, it can easily be seen that $\geq^H$ may also generate FSD, as a limiting case. This can be achieved be using the power function, giving rise to a risk-loving version of PDSD, defined as follows.

**Definition 9.** We say that $X$ dominates $Y$ w.r.t. k-PDSD and write $X \geq^{\left(1+\frac{1}{k}\right)-PD} Y$ iff $-Y \geq_{\left(1+\frac{1}{k}\right)-PD} -X$. We write $X \geq^{1-PD} Y$ iff $X \geq^{\left(1+\frac{1}{k}\right)-PD} Y$, $\forall k > 0$.

The following properties, C1–C4, follow straightforwardly.

C1) $X \geq^{\left(1+\frac{1}{k}\right)-PD} Y$ and $Y \geq_{PD}^{1+1/k} X$ iff $F_X = F_Y$.

C2) If $k_1 > k_2$, then $\geq^{\left(1+\frac{1}{k_1}\right)-PD}$ implies $\geq^{\left(1+\frac{1}{k_2}\right)-PD}$.

C3) $X \geq^{1-PD} Y$ iff $X \geq_1 Y$.

C4) $X \geq^{2-PD} Y$ iff $X \geq^2 Y$.

If $k$ is a positive integer and $U_{1:k}$ is the smallest-order statistic (or sample minimum) from a sample of i.i.d. RVs $U_1, \ldots, U_k$, the CDF of $U_{1:k}$ is $F_{U_{1:k}}(t) = 1 - (1 - F_U(t))^k$. Then, similarly to (1), condition (3) yields

$$X \geq^{\left(1+\frac{1}{k}\right)-PD} Y \text{ iff } X_{1:k} \geq^2 Y_{1:k}. \qquad (4)$$

As $k$ increases, smaller values become progressively more important, whereas the weights of larger values remain fixed.

SSD and ICX do not imply each other and can hold simultaneously, not only in the trivial case of FSD. This also holds for $\geq_{H_1}$ and $\geq^{H_2}$, for possibly different distortions $H_1, H_2$. Thus, we can introduce a definition of a mixed (risk-averse, risk-loving) DSD and extend it to the case of PDSD.

**Definition 10.**

We say that $X$ dominates $Y$ w.r.t. $H_1, H_2$-mixed-DSD and write $X \geq_{H_1}^{H_2} Y$ iff $X \geq_{H_1} Y$ and $X \geq^{H_2} Y$ hold simultaneously.

We say that $X$ dominates $Y$ w.r.t. $(1 + 1/k_1), (1 + 1/k_2)$-mixed-PDSD and write $X \geq_{1+1/k_1}^{1+1/k_2} Y$ iff $X \geq_{\left(1+\frac{1}{k_1}\right)-PD} Y$ and $X \geq^{\left(1+\frac{1}{k_2}\right)-PD} Y$ hold simultaneously.

*Example 4. symmetric case*

Let, for instance, $X$ take the values 1,2,4 and $Y$ take the values 0,2.5,3 (with uniform probabilities). It is easy to see that $X \not\geq_1 Y$. Integration yields

$$\int_0^u Q_X(p) - Q_Y(p) dp^k = \int_0^u Q_{-Y}(p) - Q_{-X}(p) dp^k = \begin{cases} u^k & 0 < u \leq \frac{1}{3} \\ \frac{3^{-k}}{2}(3 - 3^k u^k) & \frac{1}{3} < u \leq \frac{2}{3} \\ \frac{3^{-k}}{2}(3 - 3 \cdot 2^k + 2 \cdot 3^k u^k) & \frac{2}{3} < u \leq 1 \end{cases}$$

$\int_0^u Q_X(p) - Q_Y(p) dp^k = \int_0^u Q_{-Y}(p) - Q_{-X}(p) dp^k \geq 0, \forall u$ for $k \leq 1.57$ (numerical solution)

Therefore, $X \geq_{1.63}^{1.63} Y$.

## 5. Characterization of DSD through distorted expectations

Given an RV $X$ and a distortion function $\phi$, a distorted expectation, or distortion risk measure, namely $\rho_\phi(X)$, is generally defined as

$$\rho_\phi(X) = \int_{-\infty}^{\infty} t d\phi(F_X(t)) = \int_0^1 Q_X(p) d\tilde{\phi}(p) = -\int_{-\infty}^0 \tilde{\phi}(F_X(t)) dt + \int_0^{\infty} \phi(1 - F_X(t)) dt$$

where $\tilde{\phi}(t)$ is the dual distortion function of $\phi$. Clearly $\rho_\phi(X) = E(X)$ iff $\phi(t) = t$. Distorted expectations have been studied in depth in the financial and economic literature. In particular, in insurance, it is typically assumed that $X$ is a non-negative RV representing losses (Wang and Young 1998), and $\rho_\phi$ reduces to

$$\rho_\phi(X) = \int_0^{\infty} \phi(1 - F_X(t)) dt.$$

In economics, if we assume $X$ to be a non-negative RV representing income, the functional $\rho_\phi$ is particularly interesting when $\phi(t) = t^n$, where $n$ is a positive integer. In this case, we obtain the *generalized Gini indices* introduced by Donaldson and Weymark (1983), given by

$$\Xi_n(X) = \int_0^{\infty} (1 - F_X(t))^n dt, n = 1, 2, \ldots$$

where the classic Gini coefficient of inequality is given by $1 - \frac{\Xi_2(X)}{E(X)}$.

In Theorem 4, we show that $\rho_\phi$ is isotonic with FSD, SSD or ICX for every increasing, increasing convex or increasing concave distortion function $\phi$, respectively, as stated in the following theorem. Parts 1 and 3 can be obtained through Theorem 4.4 of Wang and Young (1998).

**Theorem 4**

1) $X \geq_1 Y$ iff $\rho_\phi(X) \geq \rho_\phi(Y)$ for every increasing function $\phi$.

2) $X \geq_2 Y$ iff $\rho_\phi(X) \geq \rho_\phi(Y)$ for every increasing and convex function $\phi$.

3) $X \geq^2 Y$ iff $\rho_\phi(X) \geq \rho_\phi(Y)$ for every increasing and concave function $\phi$.

Thus, all decision makers who prefer "more" to "less" may be represented, in terms of the functional $\rho_\phi$, by an increasing transformation $\phi$ (all distortion functions are increasing, by definition). Decision makers who are also risk averse may be represented by convex distortions $\phi$, that is, a restricted class of functions.

Theorem 5 characterizes $H$-DSD (as well as risk-loving $H$-DSD and mixed-DSD) in terms of isotonic distorted expectations. Intuitively, the characteristics of $\phi$ must depend on those of $H$. If $\geq_H$ is "close" to FSD, we expect that the class of distortions $\phi$ that preserve $\geq_H$ consists of "most" distortion functions, whilst, if $\geq_H$ is "close" to SSD (but stronger), the class of distortions $\phi$ that preserve $\geq_H$ consists of all convex distortions plus "some" others, which might be concave or neither convex nor concave (similar arguments hold for $\geq^H$). Generally, the weaker the order, the smaller the class. The next theorem formalizes this intuition.

**Theorem 5**

1) $X \geq_H Y$ iff $\rho_\phi(X) \geq \rho_\phi(Y)$ for every distortion function $\phi$ that is less convex than $H$ (i.e., such that $H \geq_{CX} \phi$).

2) $X \geq^H Y$ iff $\rho_\phi(X) \geq \rho_\phi(Y)$ for every distortion function $\phi$ that is more convex than $H$ (i.e., such that $\phi \geq_{CX} H$).

3) $X \geq_{H_1}^{H_2} Y$ iff $\rho_\phi(X) \geq \rho_\phi(Y)$ for every distortion function $\phi$ that is less convex than $H_1$ and more convex than $H_2$ (i.e., such that $H_1 \geq_{CX} \phi \geq_{CX} H_2$).

If $X \geq_H Y$ with $H$ convex, $\rho_\phi(X) \geq \rho_\phi(Y)$ holds for every concave function $\phi$, plus some non-concave $\phi$. For instance, let $X \geq_{1+1/k}^{PD} Y$. The distortion $\phi_m(t) = 1/2(t^{1/m} + t^m)$ is concave in the interval $[0, m^{-\frac{3m}{m^2-1}})$ and convex in $(m^{-\frac{3m}{m^2-1}}, 1]$. $\phi_m(t^k)$ is convex (i.e., $H_k \geq_{CX} \phi_m$) for $k > m$, then $X \geq_{1+1/k}^{PD} Y$ implies $\rho_{\phi_m}(X) \geq \rho_{\phi_m}(Y) \; \forall k > m$.

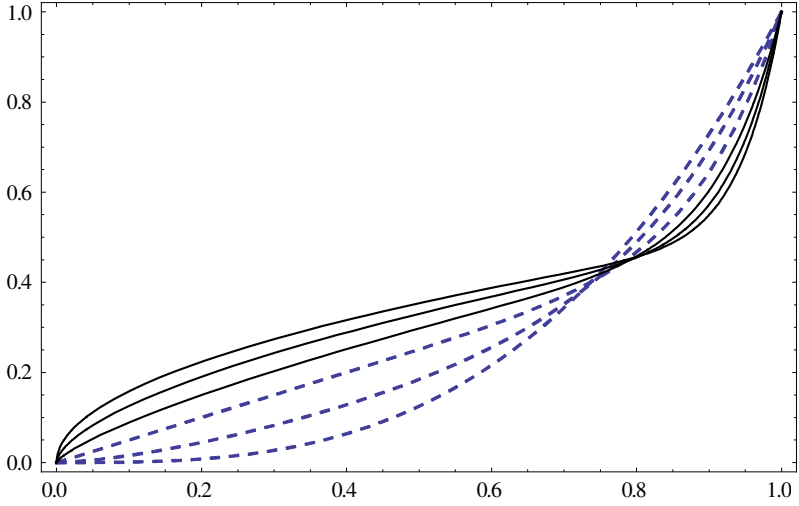

**Figure 2.** $\phi_m(t^3)$ for $m = 1,2,3$ (dashed) and for $m = 4,5,6$ (solid). $\phi_m$ is more convex than $t^3$ for $m \leq 3$.

If $H$ is concave, the class of distortions that are less convex than $H$ contains only "some" concave $\phi$s. This result conforms to the idea that $H$-DSD becomes weaker in parallel with the degree of concavity of $H$. Theorem 5 might also be stated in terms of the Arrow–Pratt measure. In fact, $\rho_\phi$ is isotonic with $H$-DSD iff $\phi$ is at most as risk averse as $H$, that is, $R_\phi(u) \leq R_H(u), \forall u \in [0,1]$. With regard to PDSD, $X \geq_{(1+\frac{1}{k})-PD} Y$ iff $\rho_\phi(X) \geq \rho_\phi(Y)$ for all distortions $\phi$ such that $R_\phi(u) \leq \frac{u}{k-1}), \forall u \in [0,1]$ (where $\frac{u}{k-1} = R_{H_k}(u)$, and $H_k(t) = t^k$).

Now, let us focus on non-negative RVs. Theorem 4 makes it possible to identify the set of generalized Gini indices $\Xi_n$ that preserve PDSD. Observe that $\Xi_n(U) = E(U_{1:n})$, where $U_{1:n} = \min\{U_1, \dots, U_n\}$ and $U_1, \dots, U_n$ are i.i.d. RVs with CDFs $F_U$. Muliere and Scarsini (1989) derived an interesting property of m-th degree inverse stochastic dominance, namely $\geq_m^{-1}$. They noted that $X \geq_m^{-1} Y$ ($m$ positive integer) implies $E(X_{1:n}) \geq E(Y_{1:n})$ for all integers $n \geq m$. Similarly, Theorem 5 states that the PDSD of order $1 + m$ ($m$ positive integer) fulfils the same property; that is, $X \geq_{(1+m)-PD} Y$ implies $E(X_{1:n}) \geq E(Y_{1:n})$ for all integers $n \geq m$.

## 6. Conclusion

We introduced a new family of stochastic orders, semiparametrized by distortion functions, generalizing FSD and SSD. The strength of the ordering relation depends on the degree of convexity of the distortion function employed. By focusing on proper classes of distortions, we may obtain a continuum of dominance relations, covering the preferences of decision makers from FSD to SSD and even beyond SSD. This can be achieved by requiring such classes to fulfil some basic conditions, described in section 2, namely C1 – antisymmetry, C2 – monotonicity, C3 – identity and C4 – consistency. We prove that the class of power distortion functions or, correspondingly, the PDSD family satisfies such properties. PDSD is particularly interesting from a statistical point of view, since it can be seen as an SSD relation between sample maxima. For single-crossing distributions $F_X, F_Y$, the PDSD of order $k$ (positive integer) is equivalent to a comparison of the expected sample maxima (of random samples of dimension $k$) from $X$ and $Y$.

We extended our approach to a risk-loving framework, enabling the generalization of ICX, yielding risk-loving-DSD and, similarly, risk-loving-PDSD. Risk-averse- and risk-loving-DSD can be combined in a "mixed" order, making it possible to represent the preferences of decision makers by controlling both their aversion and their attraction to risk.

Finally, we characterized the orders analysed in terms of distorted expectations. We derived the properties that a distortion function should fulfil to yield a distorted expectation that is isotonic with risk-averse-, risk-loving- or mixed-DSD.

Clearly, the approach used in this paper can be applied to other families of stochastic orders defined by iterated integrations, such as the Lorenz dominance of first and second degree (Aaberge 2009).

**APPENDIX**

**Proof of Proposition 1**

1) iff 2) because $Q_X(H^{-1}(t)) = Q_{X_H}$ and $Q_Y(H^{-1}(t)) = Q_{Y_H}$. 2) iff 3) by substitution. The equivalence between SSD and inverse-SSD (Thistle 1989) implies the equivalence of 3) and 4).

For the proof of Theorem 1, we need the following technical lemmas.

**Lemma A1.** Let $\varepsilon_1, \varepsilon_2, \dots, \varepsilon_n$ and $c_1, c_2, \dots, c_n$, where $c_1 \geq c_2 \geq \dots \geq c_n \geq 0$ are two sequences of real numbers. If $s_k = \sum_{j=1}^{k} \varepsilon_j \geq 0$, for every $k = 1, \dots, n$, then $\sum_{j=1}^{n} c_j \varepsilon_j \geq 0$.

*Proof*

Start with the identity $\sum_{j=1}^{n} c_j \varepsilon_j = c_1 s_1 + \dots + c_n(s_n - s_{n-1}) = s_n c_n + \sum_{j=1}^{n-1} s_j (c_j - c_{j+1})$. The thesis follows from the positivity of the terms $s_k$, $k = 1, \dots, n$ and $c_j - c_{j+1}$, $j = 1, \dots, n-1$.

**Lemma A2.** Let $r, s$ be two non-decreasing right-continuous functions on $[0,1]$ and let

$$\int_{-\infty}^{u} r(t)dt \le \int_{-\infty}^{u} s(t)dt, \forall u \in [0,1]$$

Then

$$\int_{-\infty}^{u} r(g(t))dt \le \int_{-\infty}^{u} s(g(t))dt, \forall u \in [0,1]$$

for every increasing convex function $g$.

*Proof*

By a change of variable, we need to show

$$\int_{-\infty}^{u} r(t)h(t)dt \le \int_{-\infty}^{u} s(t)h(t)dt, \forall u \in \mathbb{R}, \tag{A1}$$

for every decreasing function $h(t) > 0$. It can readily be seen that the above inequality follows from the initial assumption. Indeed, the inequality $\int_{-\infty}^{u} r(t)dt \le \int_{-\infty}^{u} s(t)dt, \forall u \in [0,1]$ implies that eq. (A1) holds for all decreasing step functions $h(t) > 0$, i.e. $h(t) = c_j$ for $a_{j-1} < t < a_j$, with a decreasing sequence of $c_j$ and an increasing sequence of $a_j, j = 1, \ldots, n, a_0 = 0, a_n = u$. This follows from Lemma A1, by setting

$$\varepsilon_j = \int_{a_{j-1}}^{a_j} (s(t) - r(t))dt,$$

because $\int_{-\infty}^{u} (s(t) - r(t))h(t)dt = \sum_{j=1}^{n} \int_{a_{j-1}}^{a_j} (s(t) - r(t))c_j dt = \sum_{j=1}^{n} c_j \varepsilon_j \ge 0$.

Finally, since a decreasing $h(t) > 0$ can be approximated by decreasing step functions, the inequality (7) holds for the general decreasing function $h(t)$.

**Proof of Theorem 1**

1) It can easily be seen that $q(u) = \frac{G'(u)}{H'(u)}$ is decreasing (i.e., $(G'/H')' \le 0$) iff $R_G(u) \ge R_H(u), \forall u \in [0,1]$. Thus, i) and ii) are equivalent. Then, the equivalence of ii) and iii) is obtained by the following relation (a similar result, stated in terms of utility functions, was proved by Pratt 1964).

$G(H^{-1}(z))'' = \left(G'(H^{-1}(z))(H^{-1})''(z) + \left((H^{-1})'(z)\right)^2 G''(H^{-1}(z))\right) \le 0$ iff

$$\frac{G''(H^{-1}(z))}{G'(H^{-1}(z))} \ge -\frac{(H^{-1})''(z)}{\left((H^{-1})'(z)\right)^2} = \frac{H''(H^{-1}(z))}{H'(H^{-1}(z))}$$

Indeed, by the inverse function theorem,

$$(H^{-1})'(z) = \frac{1}{H'(H^{-1}(z))}$$

$$(H^{-1})''(z) = \left(\frac{1}{H'(H^{-1}(z))}\right)' = -\frac{(H^{-1})'(z)H''(H^{-1}(z))}{(H'[H^{-1}(z)])^2} = -\frac{H''(H^{-1}(z))}{H'(H^{-1}(z))(H'[H^{-1}(z)])^2}$$

Thus:

$$-\frac{(H^{-1})''(z)}{\left((H^{-1})'(z)\right)^2} = \frac{H''(H^{-1}(z))}{H'(H^{-1}(z))}.$$

2) Now, we can prove that i) implies the thesis. Let $Q_X H' = M_X$, $Q_Y H' = M_Y$. We can write the condition $X \geq_H Y$ as follows:

$$\int_0^u M_X(t)dt \geq \int_0^u M_Y(t)dt, \forall u \in [0,1].$$

Then, since $q$ is decreasing, Lemma A2 yields

$$\int_0^u M_X(t)q(t)dt \geq \int_0^u M_Y(t)q(t)dt, \forall u \in [0,1].$$

**Proof of Lemma 1**

As a consequence of the existence of a finite expectation, $Q(0)H_k(0) = \lim_{t \to 0} Q(t)H_k(t) = 0$ for every $k > 1$. Thus, we obtain

$$\int_0^u Q(x)\, dH_k(x) = Q(u)H_k(u) - \int_0^u H_k(y)\, dQ(y),$$

using integration by parts.

Now let $\hat{Q}(y) = \inf\{x \in \Re |\, Q(x) \geq y\}$, the left-continuous generalized inverse of $Q$. Since $\hat{Q}(y) = F(y)$, $\lambda - a.e.$, where $\lambda$ is the Lebesgue measure, integrating by substitution (Hoffmann-Jørgensen 1994, pp.204–206), we find

$$\int_0^u H_k(x)\, dQ(x) = \int_{Q(0)}^{Q(u)} H_k\left(\hat{Q}(y)\right) dy = \int_{Q(0)}^{Q(u)} H_k(F(y))\, dy.$$

Summarizing, we have

$$\frac{1}{H_k(u)} \int_0^u Q(t)dH_k(t) = Q(u) - \int_{Q(0)}^{Q(u)} \frac{H_k(F(y))}{H_k(u)} dy.$$

For every point of continuity of $Q$, the strict inequality $H_k(F(t)) < H_k(u)$ holds whenever $t < Q(u)$. Then $\lim_{k \to \infty} \frac{H_k(F(t))}{H_k(u)} = 0$ for every $0 < t < u$. Hence, since

$$\left|\frac{H_k(F(t))}{H_k(u)}\right| = \frac{H_k(F(t))}{H_k(u)} \leq \frac{F(t)}{u}, k = 1,2,3,\ldots, 0 < t < u,$$

where $\frac{F(t)}{u}$ is an integrable function in the interval $(Q(0), Q(u))$, by the Lebesgue dominated convergence theorem, we have

$$\lim_{k \to \infty} \int_{Q(0)}^{Q(u)} \frac{H_k(F(y))}{H_k(u)} dy = \int_{Q(0)}^{Q(u)} \lim_{k \to \infty} \frac{H_k(F(t))}{H_k(u)} dy = 0.$$

**Proof of Theorem 2**

The thesis follows from Lemma 1. Clearly $X \geq_{H_k} Y$ (i.e., $\int_0^u (Q_X(t) - Q_Y(t)) dH_k(t) \geq 0, \forall u$) iff $\frac{1}{H_k(u)} \int_0^u (Q_X(t) - Q_y(t)) dH_k(t) \geq 0, \forall u$. Then, for $k \to \infty$, we obtain $X \geq_{H_k} Y$ iff $Q_X(t) - Q_Y(t) \geq 0, \forall u$.

**Proof of Theorem 5**

The thesis follows directly from Theorem 3 of Hanoch and Levi (1969).

**Proof of Theorem 4**

Point 1) can be proved easily. As for point 2), consider the following class of (increasing and) convex distortion functions:

$$g_p(z) = (p - (1 - z))_+.$$

Note that

$$\int_0^\infty g_p(1 - F_X(t)) dt$$

$$= \int_0^\infty (p - (1 - \bar{F}(t)))_+ dt = \int_0^{Q(p)} (p - F(t)) dt = pQ(p) - \int_0^{Q(p)} F(t) dt =$$

$$pQ(p) - \int_0^p z dQ(z) = pQ(p) - pQ(p) + \int_0^p Q(z) dz = \int_0^p Q(z) dz.$$

Similarly

$$\int_{-\infty}^0 \tilde{g}_p(F(t)) dt = \int_{-\infty}^0 1 - g_p(1 - F(t)) dt = \int_{-\infty}^0 1 - (p - F(t))_+ dt =$$

$$\int_{-\infty}^{Q(p)} \bigl(F(t) + (1-p)\bigr)dt - Q(p) = Q(p)(1-p) + \int_{-\infty}^{Q(p)} F(t)dt - Q(p) =$$

$$= Q(p)(1-p) + \int_0^p z\, dQ(z) - Q(p) = Q(p)(1-p) + pQ(p) - \int_0^p Q(z)dz - Q(p)$$

$$= -\int_0^p Q(z)dz.$$

Then

$$\rho_{g_p}(X) = \begin{cases} -\int_{-\infty}^0 1 - g_p(1-F(t))dt & Q(p) < 0 \\ \int_0^\infty g_p(1-F(t))dt & Q(p) > 0 \end{cases} = \int_0^p Q(z)dz.$$

Since every convex distortion function $\phi$ can be approximated by a combination of convex functions $g_p$, then $X \geq_2 Y$ iff $\rho_\phi(X) \geq \rho_\phi(Y)$, for every increasing and convex function $\phi$.

Point 3) has been proved by Wang and Young (1998) for non-negative RVs; it can easily be extended to real RVs using similar arguments to those of point 2).

**Proof of Theorem 5**

1) $X \geq_H Y$ iff $X_H \geq_2 Y_H$. Theorem 4 states that $X_H \geq_2 Y_H$ iff $\rho_\psi(X_H) \geq \rho_\psi(Y_H)$, for every increasing and convex function $\psi$. Let $\psi(1 - H(p)) = \phi(1-p)$. Then, since $1 - H^{-1}(1 - \psi^{-1}(t)) = \phi^{-1}(t)$, we obtain $\psi^{-1}(t) = 1 - H(1 - \phi^{-1}(t))$, where $\psi^{-1}$ is an increasing concave function, by construction. However, $\psi^{-1}$ is increasing concave iff $1 - \psi^{-1}(1-t) = H(\phi^{-1}(t))$ is increasing convex or, equivalently, $\phi(H^{-1}(z))$ is increasing concave. Point 2) can be proved similarly to point 1). Point 3) follows straightforwardly from points 1) and 2).